\documentclass[12pt]{article}

\usepackage{amssymb, latexsym, amsthm,amsmath,color}
\usepackage[all]{xy}
\usepackage{graphicx,psfrag,epsfig}
\usepackage[pdf]{pstricks} 

\usepackage{times}
\newtheorem{Theorem}{Theorem}[section]
\newtheorem{Proposition}[Theorem]{Proposition}
\newtheorem{Lemma}[Theorem]{Lemma}
\newtheorem{Corollary}[Theorem]{Corollary}
\newtheorem{Remark}[Theorem]{Remark}
\newtheorem{Definition}[Theorem]{Definition}

\newcommand{\bi}{\begin{enumerate}}
\newcommand{\ei}{\end{enumerate}}
\newcommand{\be}{\begin{equation}}
\newcommand{\ee}{\end{equation}}
\newcommand{\ba}{\begin{array}}
\newcommand{\ea}{\end{array}}

\def\ra{\rightarrow}

\def\t{\mathbf{t}}
\def\c{\mathbb{C}}

\def\q{\mathbb{Q}}
\def\z{\mathbb{Z}}
\def\n{\mathbb{N}}

\def\cp{{\mathcal P}}
\def\cq{{\mathcal Q}}

\usepackage{pstricks, pst-node, pst-tree} 
\input xy
\xyoption{all}

\title{Index of Kato surfaces}
\author{Akira Fujiki and Massimiliano Pontecorvo}
\date{\today}
\begin{document}

\maketitle

\begin{abstract}
The compact curves of an intermediate Kato surface $S$
form a basis of $H^2(S,\mathbb Q)$.
We present a way to compute the associated rational coefficients
of the first Chern class $c_1(S)$.
We get in particular a simple geometric obstruction for $c_1(S)$
to be an integral class, or equivalently index$(S)=1$.
We also find an expression for the exponents of the contracting germ of $S$
in terms of self-intersection numbers of the compact curves.
\end{abstract}

\section{Introduction}

A Kato surface is a minimal compact complex surface $S$ with positive
second Betti number $b_2(S)>0$ containing a global spherical shell i.e. there
is an open subset $V$ biholomorphic to a neighborhood $U$ of $S^3\subset \c^2$
with the property that $S\setminus V$ is connected.
It was shown by Kato \cite{ka77} that $\pi_1(S)\cong \z$ and that $S$ is diffeomorphic to a
Hopf surface blown up at $b_2(S)$ points. In particular $S$ admits no K\"ahler metric
and also no holomorphic pluricanonical sections
- i.e. $\mathrm{Kod}(S)=-\infty$ so that $S\in \mathrm{VII}_0^+$ in Kodaira classification.
Furthermore, Kato surfaces are the only known examples in this class and
a strong conjecture of Nakamura \cite[5.5]{na89} asserts that every $S\in \mathrm{VII}_0^+$
should be a Kato surface.
Recently, Teleman has developed instanton methods to produce curves on class-VII$^+_0$
surfaces:
it is proven in \cite{te10}\cite{te13} that every $S\in \mathrm{VII}_0^+$ with  $b_2 \leq 2$
has a cycle of rational curves and is therefore diffeomorphic to a Kato surface, by \cite[8.5]{na89};
and it is actually a Kato surface when $b_1=1$.

\medskip
As it turns out, Kato surfaces always admit exactly $b:=b_2(S)$ rational curves $D_1,...,D_b$
some of which form a cycle $C$. In the present work we will only consider
\it intermediate \rm Kato surfaces, meaning that there is at least one component
$D_i$ which is not contained in a cycle.
In this situation it is known that the maximal curve $D=\sum_{i=1}^b  D_i$
is connected and consists of a unique cycle $C$ with a number of branches
$B_1,...,B_N$ appended, $1\leq N\leq b_2(C)$.

The terminology intermediate then comes from the Dloussky number of $S$
$$\sigma(S):= - \sum_{i=1}^b D_i^2$$
which is known to satisfy the inequalities $2 b \leq \sigma(S) \leq 3b$ where
the \it extreme \rm Kato surfaces consist of the Enoki surfaces when $\sigma=2b$
and of the Inoue-Hirzebruch surfaces when $\sigma=3b$.
Both inequalities are strict precisely when $S$ has at least one branch.

We now come to the main topic of this work: apart from Enoki surfaces,
the curves $D_1,...,D_b$ form a basis for $H^2(S,\q)$ and we can write
the first Chern class of $S$, denoted by $-K \in H^2(S,\z) \setminus \{0\}$, as a linear combination
\be -K=\sum_{i=1}^b d_i D_i \quad  \mathrm{ with \;} d_i \in \q . \ee
We are interested in computing the $\mathrm{index}(S)$ which by definition is the least (positive)
integer $m$ such that $-mK$ is represented by an effective divisor;
i.e. $md_i \in\n\cup\{0\}$ for all $i$.

\medskip
Our main result is to write down explicitly all the rational coefficients $d_i$ thus solving the
linear system (1) in \cite[ p.1532]{do99} in terms of the self-intersection numbers $D_i^2$.
For doing this we make use of the Dloussky sequence $\mathrm{Dl}S$ - as described in
\cite{ot08} and \cite{dl11} - and of the dual graph of $S$ which describes the
configuration of the rational curves $D_i\subset S$.
The latter was studied by Nakamura who proved the important result that
a surface $S\in $VII$_0^+$ with $b_2(S)$ curves has the same
dual graph of some Kato surface \cite{na90}.

\medskip We concentrate on the intermediate case because
for the other Kato surfaces $-K$ is always reduced, if it can be represented by a divisor.
More precisely, Enoki surfaces satisfy $\sigma = 2b$ and
$-K$ is a divisor if and only if there is an elliptic curve,
in which case these surfaces are called parabolic Inoue and the anti-canonical  divisor
is the maximal curve which is a disjoint union of the elliptic curve with the
cycle of rational curves.
In the other case $\sigma = 3b$ we have:
hyperbolic Inoue surfaces which also have effective and
disconnected anti-canonical divisor consisting of two cycles;
the other possibility are half Inoue surfaces in which case
the first Chern class is represented by the maximal curve $C$,
just one cycle in this case, but however the anti-canonical bundle
is not a divisor.

\medskip
In section 2 we determine some fundamental properties of the multiplicities $d_i$
of $-K$ by a careful inspection of the dual graph and repeated use of adjunction formula,
starting from the case in which $S$ has only one branch.
In section 3 we explicitly compute all $d_i$'s by means of some inductively defined
multilinear forms; and in section 4 we extend the result to the general case in
which $S$ may have several branches.

In particular, our results yield topological obstructions for the existence of a numerical anti-canonical
divisor (equivalently $\mathrm{index}(S)=1$); for example we show that the number of branch components cannot exceed the number of cycle components if the index is $1$.

Finally, in section 5, we present some applications.
By a result of Apolstolov-Grantcharov-Gauduchon \cite{agg}, $\mathrm{index(S)}=1$
is a necessary condition for a Kato surface to admit a bi-Hermitian metric
and by \cite{do99} is also a necessary condition for the existence of a holomorphic vector field
$\theta\in H^0(S,\Theta_S)$. It is shown there that the zeroes of $\theta$ form
a divisor $D_{\theta}$ and our method allows to compute this divisor explicitly.
We also show how our work is related to recent results of Dloussky
who computed in \cite{dl11} the discriminant of the singularity obtained
by contracting the maximal curve of $S$ to a point.

Another important tool for the study of Kato surfaces, also introduced by Dloussky  \cite{dl84},
is a contracting holomorphic germ around the origin of $\c^2$ which gives a link
with holomorphic dynamics.
Polynomial normal forms $\phi$ for this germ were found by Favre \cite{fa00} and we show how our
method gives a way to compute the relevant exponents $j,s$ and $k$ of the polynomial $\phi$.
These exponents determine the dimension of the moduli space of logarithmic deformations of $S$ as
described in \cite{ot08} and we show in two different ways
that it equals the sum of the lengths of the regular
sequences in Dl$S$, if $S$ has no vector field.


\section{Kato surfaces with one branch}

A Kato surface $S$ with $b_2(S)=:b$ contains exactly $b$ rational curves which are
geometrically obtained as follows:
start by blowing up the origin of a ball $0\in B\subset \c^2$ at $b$ infinitely near points
in order to obtain a complex surface $\tilde B$ which will contain $b$ rational curves
of negative self-intersection the last of which $C_b$ being the only one of self-intersection $-1$.
Then, in order to get a minimal compact surface, one takes the quotient by a biholomorphism
$\psi$ which identifies a neighborhood of the origin $0$ with a small ball around some $p\in C_b$
in such a way that $p=\psi(0)$.

Notice that most curves in $S$ will have self-intersection number $-2$ because
a curve $C$ with $C^2=-(k+2)\leq -3$ can only be obtained by repeatedly blow up a
fixed node of previously created exceptional curves
and will therefore come along with a chain of $-2$-curves of length $(k-1)$.
All other curves in $S$ are either obtained in this way or else by blowing up a
general point of the previously created exceptional divisor, in particular their
self-intersection number is $-2$.
Therefore, if $(a_1,\dots ,a_b)$ denotes the string of  \it opposite \rm
self-intersection numbers of the rational curves in $S$ we can separate it
into \it singular \rm sequences $s_k:=((k+2)2\dots 2)$ of total length $k$
and \it regular \rm sequences $r_m:=(2\dots 2)$ of length $m$.

This notation was introduced by Dloussky in \cite{dl84} and the \it Dloussky sequence \rm
of a Kato surface $S$ is a sequence of $b=b_2(S)$ integers $\geq 2$
describing the opposite self-intersection number of the rational curves,
following their \it order of creation, \rm and grouped
into singular sequences and regular sequences of maximal length.
\bigskip

We consider in this section a Kato surface $S$ associated to a \it simple \rm
Dloussky sequence as described in \cite[p.335]{ot08} or \cite[p.43]{dl11}
of the following form
\be   \mathrm{DlS}=[s_{k_0}s_{k_1}...s_{k_{(p-1)}}r_m]  \label{simple} \ee
with $k_0 \geq 1$ and $m \geq 1$.
In other words $\mathrm{DlS}$ has at least one singular sequence and exactly one regular
sequence; is is equivalent to say that $S$ is of intermediate type with only one branch.
The associated dual graph $\Gamma$ represents the configuration of all curves in $S$:
each node is a rational curve with each edge connecting two nodes
whenever the corresponding curves have an intersection point.
The ($-2$)-curves will be indicated by a white node without any
further reference to their self-intersection number.
Furthermore, there are exactly $p$ curves whose self-intersection number
is not greater than $-3$; we denote them by a black node and a positive integer
indicating the opposite self-intersection number.
Our first picture illustrates this notation and points out the geometric duality between
chains of $(-2)$-curves of length $(k-1)$ and curves of self-intersection number $-(k+2)$
in a Kato surface.

\medskip
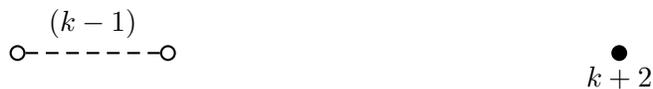
\begin{pspicture}(-6,-1)(6,1)


\cnode(-4,0){3 pt}{A}
\cnode(-2,0){3 pt}{B}


\ncline[linestyle=dashed]{-}{A}{B}
\Aput{\small{$(k-1)$}}


\qdisk(4,0){3 pt}
\uput[270](4,0){\small{$k+2$}}

\end{pspicture}
\medskip

Given a Dloussky sequence $ \mathrm{DlS}$ its dual graph is constructed by connecting
an entry with value $a$ to the entry following $a-1$ places after it on the right, in cyclic order.
When $p$ is even, the black nodes of $S$ are evenly distributed between the branch and the cycle;
when $p$ is odd, the black nodes of the cycle are one more than those in the branch.
Finally, here is a picture of the dual graph of  $\mathrm{DlS}=[s_{k_0}s_{k_1}...s_{k_{(p-1)}}r_m]$
in the case \underline{$p$ odd} and $m\geq 2$.
\bigskip

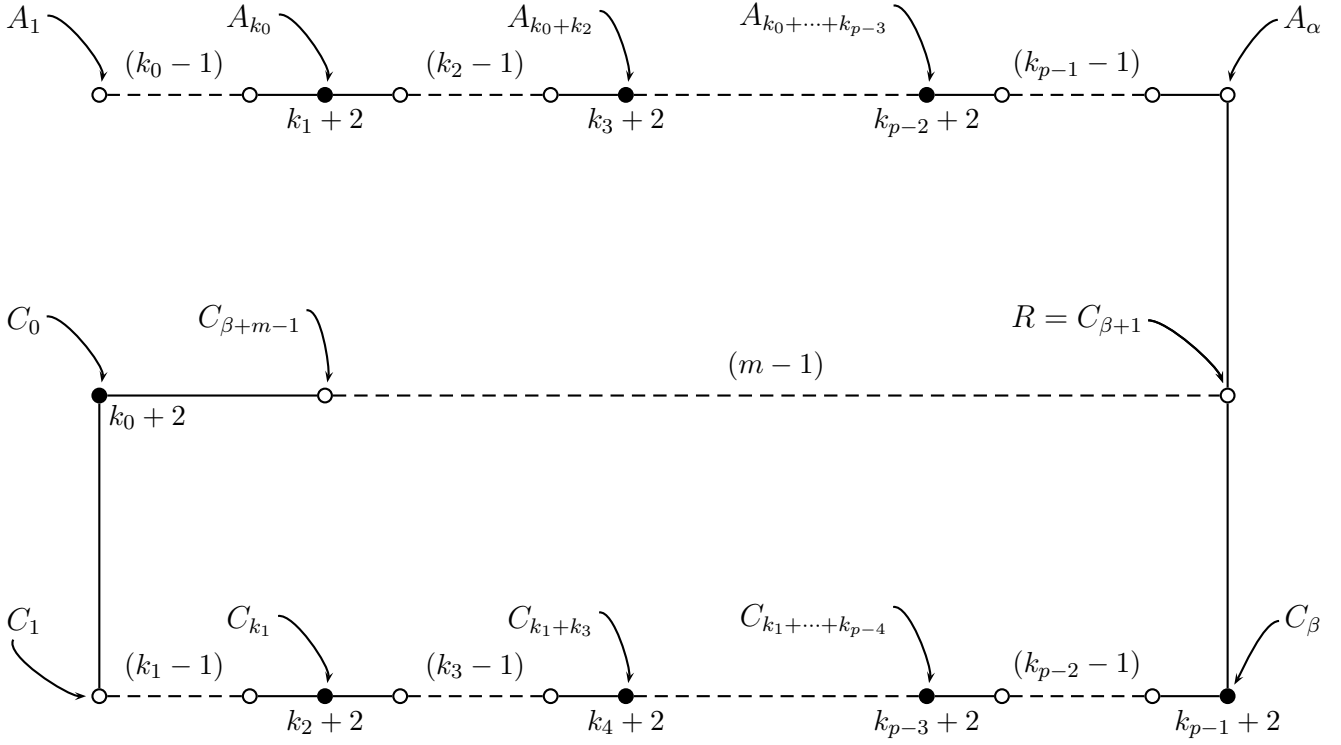
\begin{pspicture}(2,0)(20,12)

\cnode(2,10){3 pt}{A}
\rput(1,11){\rnode{a}{$A_1$}}

\cnode(4,10){3 pt}{B}

\cnode[fillstyle=solid,fillcolor=black](5,10){3 pt}{C}
\uput[d](5,10){\small{$k_1+2$}}
\rput(4,11){\rnode{c}{$A_{k_0}$}}

\cnode(6,10){3 pt}{D}

\cnode(8,10){3 pt}{E}

\cnode[fillstyle=solid,fillcolor=black](9,10){3 pt}{F}
\uput[d](9,10){\small{$k_3+2$}}
\rput(8,11){\rnode{f}{$A_{k_0+k_2}$}}

\cnode[fillstyle=solid,fillcolor=black](13,10){3 pt}{G}
\uput[d](13,10){\small{$k_{p-2}+2$}}
\rput(11.5,11){\rnode{g}{$A_{k_0+\cdots+k_{p-3}}$}}

\cnode(14,10){3 pt}{H}

\cnode(16,10){3 pt}{I}

\cnode(17,10){3 pt}{J}
\rput(18,11){\rnode{j}{$A_\alpha$}}

\ncline[linestyle=dashed]{A}{B}
\Aput{\small{$(k_0-1)$}}
\ncline{B}{C}
\ncline{C}{D}
\ncline[linestyle=dashed]{D}{E}
\Aput{\small{$(k_2-1)$}}
\ncline{E}{F}
\ncline[linestyle=dashed]{F}{G}
\ncline{G}{H}
\ncline[linestyle=dashed]{H}{I}
\Aput{\small{$(k_{p-1}-1)$}}
\ncline{I}{J}

\nccurve[ncurv=.4,angleB=10,angleA=80,nodesep=2pt]{<-}{A}{a}
\nccurve[ncurv=.4,angleB=10,angleA=80,nodesep=2pt]{<-}{C}{c}
\nccurve[ncurv=.4,angleB=10,angleA=80,nodesep=2pt]{<-}{F}{f}
\nccurve[ncurv=.4,angleB=10,angleA=80,nodesep=2pt]{<-}{G}{g}
\nccurve[ncurv=.4,angleB=180,angleA=80,nodesep=2pt]{<-}{J}{j}

\cnode[fillstyle=solid,fillcolor=black](2,6){3 pt}{A''}
\uput[dr](2,6){\small{$k_0+2$}}
\rput(1,7){\rnode{a''}{$C_0$}}

\cnode(5,6){3 pt}{B''}
\rput(4,7){\rnode{b''}{$C_{\beta+m-1}$}}

\cnode(17,6){3 pt}{C''}
\rput(15,7){\rnode{d''}{$R=C_{\beta+1}$}}
\nccurve[ncurv=.4,angleB=0,angleA=110,nodesep=2pt]{<-}{C''}{d''}

\ncline{A''}{B''}
\ncline[linestyle=dashed]{B''}{C''}
\Aput{\small{$(m-1)$}}

\nccurve[ncurv=.4,angleB=10,angleA=80,nodesep=2pt]{<-}{A''}{a''}
\nccurve[ncurv=.4,angleB=10,angleA=80,nodesep=2pt]{<-}{B''}{b''}
\nccurve[ncurv=.4,angleB=0,angleA=110,nodesep=2pt]{<-}{C''}{d''}

\cnode(2,2){3 pt}{A'}
\rput(1,3){\rnode{a'}{$C_1$}}

\cnode(4,2){3 pt}{B'}

\cnode[fillstyle=solid,fillcolor=black](5,2){3 pt}{C'}
\uput[d](5,2){\small{$k_2+2$}}
\rput(4,3){\rnode{c'}{$C_{k_1}$}}

\cnode(6,2){3 pt}{D'}

\cnode(8,2){3 pt}{E'}

\cnode[fillstyle=solid,fillcolor=black](9,2){3 pt}{F'}
\uput[d](9,2){\small{$k_4+2$}}
\rput(8,3){\rnode{f'}{$C_{k_1+k_3}$}}

\cnode[fillstyle=solid,fillcolor=black](13,2){3 pt}{G'}
\uput[d](13,2){\small{$k_{p-3}+2$}}
\rput(11.5,3){\rnode{g'}{$C_{k_1+\cdots +k_{p-4}}$}}

\cnode(14,2){3 pt}{H'}

\cnode(16,2){3 pt}{I'}

\cnode[fillstyle=solid,fillcolor=black](17,2){3 pt}{J'}
\uput[d](17,2){\small{$k_{p-1}+2$}}
\rput(18,3){\rnode{j'}{$C_{\beta}$}}

\ncline[linestyle=dashed]{A'}{B'}
\Aput{\small{$(k_1-1)$}}
\ncline{B'}{C'}
\ncline{C'}{D'}
\ncline[linestyle=dashed]{D'}{E'}
\Aput{\small{$(k_3-1)$}}
\ncline{E'}{F'}
\ncline[linestyle=dashed]{F'}{G'}
\ncline{G'}{H'}
\ncline[linestyle=dashed]{H'}{I'}
\Aput{\small{$(k_{p-2}-1)$}}
\ncline{I'}{J'}

\ncline{A'}{A''}
\ncline{J}{C''}
\ncline{J'}{C''}

\nccurve[ncurv=.4,angleB=240,angleA=180,nodesep=2pt]{<-}{A'}{a'}
\nccurve[ncurv=.4,angleB=10,angleA=80,nodesep=2pt]{<-}{C'}{c'}
\nccurve[ncurv=.4,angleB=10,angleA=80,nodesep=2pt]{<-}{F'}{f'}
\nccurve[ncurv=.4,angleB=10,angleA=80,nodesep=2pt]{<-}{G'}{g'}
\nccurve[ncurv=.4,angleB=180,angleA=60,nodesep=2pt]{<-}{J'}{j'}

\end{pspicture}

\bigskip

$\Gamma$ is usually called a \it directed \rm dual graph, meaning that it is weighted
by the self-intersection numbers and that the cycle has a given orientation.
The notations are as follows: if $D$ denotes the maximal curve in $S$
we have that $D$ has $b=b_2(S)$ irreducible components and set
$$\alpha=\sum_{0\leq 2i <p } k_{2i} \textrm{\hspace{2cm} and \hspace{2cm} }
 \beta=\sum_{1\leq 2i+1<p} k_{(2i+1)} .$$
Then $S$ has a unique branch $B$ with $b_2(B)=\alpha$ and
the irreducible components of $B$ are ordered from the tip to the root
and excluding it:
$A_1,...,A_{k_0},...,A_{k_0 +k_2},...,A_\alpha.$
Their self-intersection numbers are, respectively:
$-2,...,-(k_1 +2),...,-(k_3 +2),...,-(k_{p-2}+2)$
in the $p$ odd case of the picture.

For the unique cycle $C$ of $S$ we have: $b_2(C)=\beta +m$
and we ordered the irreducible components cyclically starting from the
one of self-intersection number $-(k_0 +2)$ which we denote by $C_0$
then going counterclockwise we encounter all other black nodes
$C_{k_1},C_{k_1 +k_3},...,C_{\beta}$ with self-intersection numbers
$-(k_2 +2),-(k_4 +2),...,-(k_{p-1} +2)$ respectively, until we reach the root $R:=C_{\beta+1}$;
after that comes a chain of $(-2)$--curves of length $m-1$ and we are back
at the black node $C_0$, if $m\geq 2$.
The case $m=1$ is special because the root is a black node
coinciding with $C_0$: $R=C_{\beta +1}=C_0$.

Also notice that in general,
$b:=b_2(S)=b_2(B)+b_2(C)=\alpha + \beta + 1 + m -1 =\alpha + \beta + m$.

\bigskip
We are interested in the first Chern class of $S$ which we denote by
$-K\in H^2(S,\z)$. Because the irreducible components of $D$ form
a basis of $H^2(S,\q)$ we know a priori that there exist rational
coefficients which we call \it multiplicities \rm
$a_1,...,a_\alpha,c_0,...,c_\beta,...,c_{\beta+m-1}$
such that the following equation holds in cohomology
$$-K=\sum_{i=1}^{\alpha} a_i A_i + \sum_{j=0}^{\beta +m-1} c_j C_j $$

This equation is completely equivalent to the linear system
\cite[(1) p.1532]{do99} via the adjunction formula applied to every
irreducible component $D_i$ of the maximal curve $D$:
\be 2=(-K-D_i)D_i. \ee

\bigskip
Suppose that $(D -D_i)D_i =2$, in other words $D_i$ is any irreducible
component different from the root $C_{\beta +1}$ and different from the
tip $A_1$. Let $D_{i-1}$ and $D_{i+1}$ be the components immediately preceding
and immediately following $D_i$ in the order illustrated by the picture.
Let finally $d_{i-1}, d_i, d_{i+1}$ denote their multiplicities.
The following two lemmas are a straightforward consequence of adjunction
and will prove to be very useful.

\begin{Lemma}{\bf (White lemma.)} Let $D_i$ be a white node ($D_i^2=-2$)
different from a tip or a root, then adjunction formula reads
$$d_{i+1} - d_i = d_i - d_{i-1}$$
In other words the multiplicities of a chain of $(-2)$-curves grow linearly.
\end{Lemma}

\begin{Lemma}{\bf (Black lemma.)} If $D_i$ is a black node with $D_i^2=-(k+2)$
and $(D -D_i)D_i =2$, the adjunction formula reads
$$d_{i+1} - d_i = d_i - d_{i-1} +k(d_i -1)$$
\end{Lemma}

\begin{Remark}\label{bw}  \rm
We are interested in measuring the growth of the coefficients of $-K$ and we just
observed that the difference $d_i -d_{i-1}$ remains constant at white nodes so
that the multiplicity $d_i$ is a piecewise linear function of ${i}$ whose slope
changes precisely at every black node.
We will use the following notation for these slopes:
 $g_0:=a_2-a_1$ ; $g_2:=a_{{k_{0}+1}}-a_{k_0}$ and more generally at
 the branch black node $A_{k_0+...+k_{2j}}$ we set:
 $ g_{2(j+1)} := a_{k_0+...+k_{2j} +1} - a_{k_0+...+k_{2j}} $ and rewrite
  Black lemma as:
\be g_{2(j+1)} - g_{2j} = k_{2j+1}(a_{k_0+...+k_{2j}}-1 ).\ee
For the cycle black nodes we set:
$g_1:=c_1-c_0$ and $g_{2j+1}:=c_{k_1+...+k_{2j+1}+1}-c_{k_1+...+k_{2j+1}}$
and rewrite Black lemma at the cycle black node $C_{k_1+...+k_{2j+1}} \neq C_0$ as
\be g_{2j+1} - g_{2j-1} = k_{2j}(c_{k_1+...+k_{2j+1}}-1 ).\ee

Finally, at the black node $C_0$ we will show in the proof of next theorem that for all $m \geq 1 $
\be  g_1 +1 = k_0(c_0 -1).  \ee
\hfill $\bigtriangleup$    \end{Remark}

From now on we denote by $R$ the root $C_{\beta +1}$ and use $r$ for its multiplicity $c_{\beta +1}$.
As an application of White and Black lemma we have:

\begin{Proposition} For all $m\geq 1$, the root multiplicity $r:=c_{\beta +1}$ and all
the multiplicities $a_i$ of the branch $B$, for $2\leq i \leq \alpha$,
can be computed in terms  of the multiplicity $a_1$ of the tip.
\end{Proposition}
\it Proof. \rm We set $\t:=a_1$ and prove that all other coefficients can be expressed
as a (piecewise linear) function of $\t$.
We start by applying adjunction formula to the tip $A_1$, assuming for simplicity that $A_1^2=-2$,
or equivalently $k_0>1$: \newline
$2=(-K-A_1)A_1=(\t-1)A_1^2 +a_2=-2\t +2 +a_2$ so that $a_2=2\t$ and $a_2-a_1=\t$. 
Now we show how these ``initial conditions" $\t:=a_1=a_2-a_1$ uniquely determine
all other coefficients. Setting $g_0:=a_2-a_1$ we have by White lemma that
$a_i=g_0\cdot i$ for $i=1,...,k_0$. We then apply Black lemma to the black
node $A_{k_0}$ so that $g_2:=a_{k_0 +1} -a_{k_0}=g_0 +k_1(a_{k_0} -1)$ which we
can explicitly compute from $g_0=\t$ and $a_{k_0}=k_0 \t$. Then by White lemma
$a_{k_0 +i}=a_{k_0} +g_2\cdot i$ for all $i=1,...,k_2$. Proceeding in this way,
the White and Black lemmas produce all the multiplicities up to $a_\alpha$
and in fact up to the root multiplicity $r=c_{\beta +1}$,
as (piecewise) linear functions of $\t$.
\hfill $\Box$

\medskip
A similar, more involved argument, produces all other multiplicities as well.

\begin{Theorem} \label{chase}
The multiplicities of the cycle $C$ of a Kato surface $S$ with
$ \mathrm{DlS}=[s_{k_0}s_{k_1}...s_{k_{p-1}}r_m] $,
are uniquely determined by the following initial conditions:
$c_0=\t+1$ and $c_1 -c_0 = k_0 \t -1$.
Where, $\t:=a_1$ is the multiplicity of the tip.
Furthermore, when $m\geq 2$, we have $r-c_{\beta +1+i}=i$ for all $i=0,1,\cdots, m-2$.

\end{Theorem}
\it Proof. \rm If we set $c_0=\t+1$ and $g_1:=c_1 -c_0 =k_0 \t -1$,
by White lemma we get all multiplicities up to $c_{k_1}$ as follows:
$c_i = c_0 + g_1\cdot i $ for each $i=1,...,k_1$.
We then apply Black lemma at the black node $C_{k_1}$ and get the
slope of the next line on which the multiplicities lie: setting
$g_3:= c_{k_1 +1} -c_{k_1}$, by Black lemma we have $g_3=g_1-k_2(c_{k_1} -1)$.
White lemma then yields $c_{k_1 +i} = c_{k_1} + g_3 i $ for each $i=1,...,k_3$.
Continuing this way, we get the cycle multiplicities $c_i$ for $i\leq\beta +1$,
all the way up to the root; this procedure is exactly the same procedure we applied
to the branch, it implies the following important lemma which we
will need for finishing the proof and is of independent interest because
it shows how the geometric
duality between self-intersection numbers of black nodes in the branch and lengths
of chains of $(-2)$-curves in the cycle (and vice versa) is reflected
in an arithmetic duality between multiplicities $a_{k_0}, a_{k_0 +k_2},...,a_{\alpha}$
of black nodes in the branch and slopes $g_1,g_3,...,g_{\beta}$ on the cycle; and vice versa.

\begin{Lemma} \label{dl} {\bf (Duality lemma.)} Consider a simple Dloussky sequence
$ \mathrm{DlS}=[s_{k_0}s_{k_1}...s_{k_{p-1}}r_m] $ with $m\geq 1$.
The above initial conditions:
$c_0=\t+1$ and $g_1=k_0 \t -1$ where $\t=a_1=g_0$ produce
 the following relations among slopes in the branch $B$
and multiplicities in the cycle $C$, for every $j$ such that $0\leq 2j < p$
$$ 1+ g_{2j} = c_{k_1+k_3+\cdots +k_{2j -1}} $$
and vice versa:
$$ 1+ g_{2j+1} = a_{k_0+k_2+\cdots +k_{2j}}. $$
\end{Lemma}

Using this lemma, which we will prove after finishing the proof of theorem,
we can now show the following:

(i) the root multiplicity $c_{\beta +1}$ obtained by adjunction applied to the cycle node
    $C_{\beta}$ coincides with the multiplicity produced by adjunction applied to
    the last branch node $A_\alpha$.

(ii) compute the multiplicities of the last chain of $(-2)$-curves:  
      $c_{\beta +i}$, for $2\leq i \leq (m-1)$ and $m\geq 2$.

(iii) check that adjunction formula holds at the root $C_{\beta +1}$
      and at the first black node $C_0$.
\medskip

In fact, (i) and (iii) assure us that we have found the unique solution, in terms of $\t$.

Starting from the proof of (i), we distinguish two cases. When $p=2q+1$ is odd we have
$\alpha=2q$ and $\beta=2q-1$. Because $A_\alpha$ is a white node, when we compute
the root multiplicity $c_{\beta +1}=r$ by adjunction formula at this last branch node
we get $r=a_\alpha +g_{2q}$; because $C_\beta$ is a black node
adjunction applied to the cycle will give $r=c_\beta +g_{2q+1}$.
Therefore (i) holds because $g_{2q+1}=a_\alpha -1$ and $c_\beta = g_{2q} +1$ by duality.

In the other case $p=2q$ is even; $\alpha=2q-2$ and $\beta=2q-1$.
Now, $A_\alpha$ is black and $r=a_\alpha +g_{2q}$ from the branch.
While $C_\beta$ is white and therefore $r=c_\beta+ g_{2q-1}$ from the cycle.
These two values agree because $g_{2q-1} = a_\alpha -1$ and $c_\beta =g_{2q}+1$
by duality.

This proves (i) and notice that we also have the following useful identity:
\be r = a_\alpha + c_\beta -1,
    \quad\mathrm{ for \; all \;} p \quad\mathrm { and \; all \;} m \geq 1.  \label{star} \ee
In fact, for every $p$ even or odd,
$r-a_\alpha -c_\beta = g_{2q} -c_\beta = g_{2q} - c_{k_1+\cdots +k_{2q-1}}$
which by duality equals $-1$.

\smallskip
(ii) We now compute the multiplicities in the last chain:
 $C_{\beta +1}=R, C_{\beta +2},\dots ,C_{\beta +m-1}, C_{\beta +m}=C_0$.

When $m\geq 2$, we can compute $c_{\beta +2}$ by adjunction to the root $R=C_{\beta +1}$
    which is a white node:       
   $2=(-K-R)R=a_{\alpha} +c_{\beta} +c_{\beta +2} -2(-1 +r)$
   from which we get that        
   $c_{\beta +2} -r = r -c_{\beta}  -a_{\alpha} = -1$.

Therefore we always have the following remarkable identity which does not depend on $\t$
\be \label{remarkable}
c_{\beta +2} = -1 +c_{\beta +1}
\ee

By White lemma, it then follows that for $m\geq 2$ and each $1\leq i \leq m -1$ we have
\be
c_{\beta + 1+i} = -i +c_{\beta +1}
\ee
in particular,
\be
c_{\beta +m-1} = c_0 +1 = \t+2
\ee

(iii) Adjunction formula at the root $R$ has been already used several times in our
construction and  therefore holds automatically.
Therefore, it only remains to check adjunction formula at the node $C_0$
(equivalently, Black lemma holds).
We distinguish two cases, suppose first that $m\geq 2$
so that $C_0\neq R$ meets only two other components;
adjunction then reads
$2=(-K-C_0)C_0=c_1 +c_{\beta +m-1} -(k_0+2)(c_0-1)$ giving
$c_1-c_0=c_0 -c_{\beta -m-1} +k_0(c_0 -1)$
which is equivalent to $g_1=-1 +k_0(c_0 -1)$ because of (\ref{remarkable}).
Finally, from the initial condition $c_0 -1 =\t$ we get $g_1=1- +a_{k_{0}}$
which is identically true.

It remains to see the case $m=1$ in which $C_0=R$ is the root
and therefore meets three irreducible components. By adjunction:
$2=(-K-C_0)C_0=c_\beta +a_\alpha +c_1 -(k_0+2)(c_0 -1)$ therefore
$0=c_1-c_0+c_\beta +a_\alpha -c_0 -k_0(c_0 -1)
=g_1 +1 -k_0(c_0 -1)$ by (\ref{star}); this is also
an identity for every $\t$ because of the initial conditions $g_1=k_0\t-1$ and $c_0=\t +1$.
\hfill $\Box$  \bigskip

\it Proof of Lemma 2.6. \rm With the given initial conditions,
 Black lemma at $A_{k_0}$ yields
 $g_2 = g_0 +k_1(k_0 a_1 -1)= c_0 -1 +k_1 g_1 = c_{k_1} -1$.
 Black lemma at $C_{k_1}$ reads:
 $g_3 = g_1 +k_2(c_{k_1} -1)= k_0 \t -1 + k_2 g_2 = -1 +a_{k_0 + k_2}$
 so that both formulas hold for $j=1$; and also for $j=0$ if we set $k_{-1}=0$.
 Supposing that they hold at $j-1$, we prove them at $j$:
 by  Black lemma at $A_{k_0 +\cdots +k_{2j-2}}$ we have
 $g_{2j} =g_{2j-2} + k_{2j-1}(-1 +a_{k_0 +\cdots +k_{2j-2}})$
which by induction equals to $g_{2j-2} + k_{2j-1} g_{2j-1} =
 g_{2j-2} + c_{k_1+k_3+\cdots +k_{2j -1}} - c_{k_1+k_3+\cdots +k_{2j -3}}=
 -1 + c_{k_1+k_3+\cdots +k_{2j -1}}$ again, by induction.
 Therefore we have proved the first identity.

 Black lemma at $C_{k_1+k_3+\cdots +k_{(2j -1)}}$ reads:
 $g_{2j +1} =g_{2j-1} + k_{2j}(-1 +c_{k_1+\cdots +k_{2j -1}})$
which by what we just proved and induction equals
 $a_{k_0 +\cdots +k_{2j-2}} -1 +k_{2j} g_{2j}=-1 +a_{k_0 +\cdots +k_{2j}}$
 \hfill $\Box$ \medskip

\bigskip
The Duality lemma tells us that on a Kato surface with one branch the branch multiplicities
determine the cycle multiplicities;
coupled with Remark \ref{bw} it also immediately yields the following useful formula.

\begin{Corollary} \label{gee} On an intermediate Kato surface with
one branch the following holds.   For all $1\leq j \leq p$
\be g_{j+1} - g_{j-1} = k_j g_j \ee
\end{Corollary}

We also point out the following identity which follows immediately from (\ref{star})
and duality.

\begin{Corollary} \label{root}
The multiplicity of the root of a Kato surface $S$ with
$\mathrm{DlS}=[s_{k_0}s_{k_1}...s_{k_{p-1}}r_m] $ and $m\geq 1$ satisfies,
$$r = g_p +g_{p-1} +1$$
\end{Corollary}


\section{Index}

In order to conveniently write in closed form all the coefficients
of the first Chern class $-K$ of an intermediate Kato surface $S$,
we now introduce some multilinear forms.
For simplicity of exposition, we will assume in this section that $S$
has only one branch.

It is clear from Duality lemma that we only need to find an
expression for the slopes $g_j$'s which, by Corollary \ref{gee}
can be computed inductively.
Once we have expressed $g_j$ for $ j\in\{0,...,p\}$
we can also compute the root multiplicity by means of Corollary \ref{root}.
This plan suggests the following

\begin{Definition} \rm Let $X_1, X_2,...,X_n$ denote a set of variables
and define polynomials $f$ in $n$ variables inductively, by
$$f(X_1):= X_1,   \quad f(X_1,X_2):=X_1X_2 +1 $$
\be f(X_1,...,X_n):=X_n f(X_1,...,X_{n-1}) + f(X_1,...,X_{n-2}).
\label{f} \ee
 We also introduce multilinear forms $\cp$, inductively defined from $f$ as follows
 $$  \mathcal{P}(X_1) = X_1, \quad  \mathcal{P}(X_1,X_2)=X_1X_2+X_1 $$
\be \cp(X_1,...,X_n) := X_n  f(X_1,...,X_{n-1}) + \cp(X_1,...,X_{n-1})
\label{p}\ee
 \hfill     $\bigtriangleup$  \end{Definition}

Notice that
$$\mathcal{P} (X_1,...,X_n) =
X_1\cdots X_n + X_1\cdots X_{n-1} + \mbox{\small{ lower order terms. }}$$
furthermore, of course: $\mathcal{P} (X_1,...,X_{n-1}) = \mathcal{P} (X_1,...,X_{n-1},0)$.

\medskip
In order to make the notation clearer we will now use $r(p)$ for the multiplicity
of the root of a simple Dloussky sequence with $p$ black nodes
which was previously denoted by $r$ or $c_{\beta +1}$.

\begin{Proposition} \label{explicit}
Let $\mathrm{DlS}=[s_{k_0}s_{k_1}...s_{k_{p-1}}r_m]$ be a simple Dloussky sequence
with $m\geq 1$. As usual, $\t:=a_1$ denotes the multiplicity of the tip.
Then, for the slopes $g_j$'s and for the multiplicity
$r(p)$ of the root we have the following formulas which
express them as linear function of $\t$ with coefficients depending on the
black nodes only:

$$g_j = f(k_0,...,k_{j-1}) \t - f(k_1,...,k_{j-1}), \quad \mathrm{ for \; all } \quad  j = 0,1,2,...,p.$$

$$r(p) = [ \cp(k_0,k_1,...,k_{p-1}) +1] \t - \cp(k_1,k_2,...,k_{p-1}).$$
\end{Proposition}

\it Proof. \rm The proof is by induction on the number $p$ of black nodes.
Recall from Theorem \ref{chase} that $g_0=\t$ and $g_1=k_0\t-1$,
so that the definition of $f_0$ and $f(X_1)$
is consistent  with our initial conditions.
Finally, $f$ is defined by induction in such a way that $g_j$ will satisfy
Corollary \ref{gee} automatically.

In a similar way we can verify the formula for $r(p)$.
The first step of induction holds because from the final part of the proof of Theorem \ref{chase}
  $$r(1)=c_0 +g_1=a_{k_0}+g_0=(k_0 +1)\t = [\cp(k_0) +1]\t.$$

From Corollary \ref{root} we have $r(p)=g_p - g_{p-1} +1$.
It then follows from Corollary \ref{gee} that
$r(p+1)-r(p) = g_{p+1} - g_{p-1} = k_p g_p$
and to complete the proof it is enough to recall the definition of $\cp$.
\hfill $\Box$ \bigskip

\begin{Remark} \rm
Notice that in the above formulas the coefficient of $\t$ depends
on all the black nodes, while the constant coefficient does not depend
on the first black node.
  \hfill $\bigtriangleup$
\end{Remark}
\medskip

Now that we know how to express the root multiplicity  -
and in fact all multiplicities - in terms of the tip multiplicity,
we can show how to compute $a_1=\t$ and find the index of $S$;
we start by recalling this important notion.

\begin{Definition} \rm   The index of a Kato surface $S$ is the smallest
integer $m$ such that the cohomology class $-m K \in H^2(S,\z)$
is represented by an effective divisor.
In what follows $\mathrm{index}(S)$ shall denote the index of the Kato surface $S$.
\hfill $\bigtriangleup$
\end{Definition}

Therefore $\mathrm{index}(S)$ is the least common multiple of all denominators
of the rational coefficients of the first Chern class $-K=c_1(S)$.
 Because we have just expressed all multiplicities as functions of $\t$,
 we now have:

\begin{Theorem} \label{t}
Let $\mathrm{DlS}=[s_{k_0}s_{k_1}...s_{k_{p-1}}r_m]$ be a simple Dloussky sequence
with $m\geq 1$.
 The rational number $\t:=a_1$ is the solution of the linear equation  $\t +m = r(p)$.
We can then write all the multiplicities of the first Chern class $-K$ explicitly,
by substituting:
$$\t= \frac{m+\cp(k_1,k_2,...,k_{p-1})}{\cp(k_0,k_1,...,k_{p-1})}. $$
In particular, we always have $\t >0$ and furthermore every multiplicity of $-K$
is a positive integer if and only if $\t\in\n$.
Finally,
$$\mathrm{index}(S)= \frac{\cp(k_0,k_1,...,k_{p-1}) }
  {\mathrm{g.c.d}\{ m+\cp(k_1,k_2,...,k_{p-1}) , \cp(k_0,k_1,...,k_{p-1})\}}$$
\end{Theorem}
\it Proof. \rm 
The coefficients of $-K$ that we wrote down satisfy adjunction formula
at every node, for every value of $\t:= a_1$.
In order to find $\t$ we simply remark that the value of the root multiplicity
$r(p)$ (or $c_{\beta +1}$) which was found before -- starting from the tip $a_1$
and going up on the branch; or else starting from the first black node $C_0$
in the cycle and going counterclockwise --
has to coincide with the value which we get starting from $C_0$ and
going clockwise around the cycle until we hit the root $C_{\beta +1}$.

This computation is a lot easier because $C_0, C_{\beta +m-1},...,C_{\beta +1}$
is a chain of $(-2)$--curves.
During the proof of Theorem \ref{chase} we have shown the remarkable result that
for $m\geq 2$
\be
 c_{\beta +m-1} -c_0 =1.
\ee
Therefore, the fact that $c_0 =\t +1$, which holds for all $m \geq 1$,
 together with White lemma immediately give that
$$r(p) =c_{\beta +1} = \t+m$$
as wanted. Notice that this holds for $m=1$ as well -- i.e.
the root is a black node: $C_0=C_{\beta +1}$.

To see the statement about the index of $S$, just recall that $\t=a_1$
and that, using Proposition \ref{explicit} and Duality lemma,
we can write down explicitly all other coefficients of $-K$ as linear forms
in the variable $\t$ with integer coefficients,
which are polynomials in $k_0,...,k_{p-1}$ and $m$.
\hfill $\Box$

\begin{Remark}  \rm
When $p=1$ we have $\mathrm{DlS}=[s_{k_0} r_m]$ and the root  multiplicity is
$r(1)= (k_0 +1)\t$.     Therefore $\t$ solves the equation
$$\t+m =(k_0 +1)\t \quad \Leftrightarrow \quad \t=m/k_0$$
Because in this case $b_2(B)=k_0$ and $b_2(C)=m$,
we get that $\mathrm{index}(S)=1$ if and only if the number of irreducible components
of the branch divides the number of irreducible components of the cycle.

When $p=2$, it easily turns out that
$$\t=(m +k_1)/k_0(k_1 +1).$$
In this case $b_2(B)=k_0$ and $b_2(C)=k_1 +m$ so that both $k_0$ and $k_1 +1$
divide the number of cycle components if the index is $1$.

When $p=3$ the index of $S$ is the denominator of the rational number
$$ (k_2k_1 +k_1 +m )/( k_2k_1k_0 +k_1k_0 +k_2 +k_0) $$
 \hfill     $\bigtriangleup$
 \end{Remark}  \smallskip

Because for $p\geq 3$ the multilinear forms $\cp(k_0,...k_{p-1})$ are irreducible \cite{dl11};
it becomes increasingly difficult to draw precise geometrical consequences from Theorem
\ref{t}; in what follows we present a few necessary conditions for an intermediate Kato
surface to have index $1$, which however are far from being sufficient.

The following result is important for understanding the behavior of the multiplicities
of $-K$ when the index is $1$.

\begin{Corollary}\label{slopes}
Let $S$ be an intermediate Kato surface with
$\mathrm{DlS}=[s_{k_0}s_{k_1}...s_{k_{p-1}}r_m]$ and suppose that
$\mathrm{index}(S)=1$ then, the slopes $g_j \geq 1$ for all $j$ with the only
possible exception $g_1=0$, which occurs precisely for $\t =k_0=1$.
In particular, in the directed dual graph $\Gamma$, the multiplicities of $-K$
decrease only along the last chain of $(m-1)$ white nodes
going from the root $R$ back to the first black node $C_0$; see (\ref{remarkable}) .
\end{Corollary}
\it Proof. \rm First of all, $\t \in \n$ implies $g_j \geq 0$ for all $j$ by Corollary \ref{gee}
because $g_0=\t$ and $g_1=k_0\t -1 \geq 0$.
In particular, $g_{2j}$ is always strictly positive as well as  $g_{2j+1}$ for $j\geq 1$.
Finally, $g_1=k_0\t -1$ vanishes if and only if  $\t =k_0=1$.
\hfill    $\Box$
\medskip

As a first consequence we have:

\begin{Proposition} \label{broot}
Let $S$ be a Kato surface with one branch and $m=1$.
Then $\mathrm{index}(S)=1$ if and only if
$$\mathrm{DlS}=[3 s_{k_1} 2]  \quad for \; some \;  k_1\geq 0.$$
In particular, $R^2=-3$ and the branch has a unique irreducible component $A$;
furthermore $-K=A+2C$ where $C$ denotes the cycle of $S$ and its Dloussky
number satisfies $\sigma(S) = 3b_2(S) -1$.
\end{Proposition}
\it Proof. \rm This is the case in which the root is a black node
or equivalently the last chain of $(m-1)$ white nodes is empty.
By Corollary \ref{slopes} the cycle multiplicities are always non-decreasing
and therefore must be constant.
It follows that $g_{2j+1}=0$ for all $j\geq 0$, forcing
$\t =k_0=1$ and $k_{2}=0$ by Corollary \ref{gee}.
The last assertions are easily verified because
$r=c_0=\t +1=2$ by Theorem \ref{chase}
\hfill    $\Box$
\medskip

We also have the following geometric application which in the next section
we will show to hold for any intermediate Kato surface.

\begin{Proposition}\label{shortbranch}
 Let $S$ be a Kato surface with one branch,
$\mathrm{DlS}=[s_{k_0}s_{k_1}...s_{k_{p-1}}r_m]$ and suppose
$\mathrm{index}(S)=1$; then $b_2(B)\leq m$ with equality if and only if
$\mathrm{DlS}=[s_m r_m]$ or $\mathrm{DlS}=[3 s_{k_1} 2]$.
In particular, the number of branch components
cannot exceed the number of cycle components when the index is $1$.
\end{Proposition}
\it Proof. \rm
The statement certainly holds when $\mathrm{DlS}=[3 s_{k_1} 2]$
because in this case $b_2(B)=1=m$.
We can therefore assume $m\geq 2$. Let us recall that $a_1=\t$, $c_0=\t+1$ and that
$c_{\beta +m-1} -c_0 =...=r-c_{\beta +2} = 1$.
Therefore, the sequence of integers $\t,c_0,c_{\beta +m-1},...,c_{\beta +2},r$
has $m+1$ elements which are increasing - as slow as possible, $1$ by $1$ -
from the value $\t$ up to $r$.

We compare it with another strictly increasing sequence of integers
$a_1,a_2,...,a_{\alpha}, r$: which has length $\alpha +1=b_2(B) +1$ and
the same end points but grows piecewise linearly with slopes  $g_{2j} \geq 1$,
according to Corollary \ref{slopes}.
It follows that $b_2(B) +1 \leq m+1$ with equality if and only if $g_{2j}=1$.
This can only happen for $j=0$ and $\t=1$ or $j=1$ and we get from Corollary \ref{gee}
either $k_1=0$ or $g_1=0$. In the first case Dl$S=[s_{k_0}r_m]$ with $1=\t=\frac{m}{k_0}$
which is the case Dl$S=[s_mr_m]$; otherwise $k_0=1$ and $1=\t=\frac{m +k_1}{k_1 +1}$
so that $m=1$ and the proof is complete.
\hfill    $\Box$
\medskip

\begin{Remark}  \rm
 Let $S$ be a Kato surface with one branch and suppose
$\mathrm{index}(S)=1$. This happens if and only if $\t$ is an integer which is
automatically positive by Theorem \ref{t}.
We have also shown that all other multiplicities are strictly bigger than $\t$,
and therefore $\geq 2$,
except for the tip $a_1 =\t$ (which can possibly be $1$)
and is always the strictly minimal multiplicity in $S$.
The maximal multiplicity occurs at the root $R$,
while the minimal multiplicity occurs at the first cycle black node $C_0$ where $c_0=\t+1$;
both multiplicities are strict extrema in the cycle if and only if $k_0\t >1$.
For example, a black root is not a strict maximum if the index is $1$
and the surface has only one branch.

In particular $-K$ is always represented by a highly non-reduced divisor which is
strictly bigger than the maximal curve $D$, this was already known \cite{dl06} \cite{do99}.
\hfill $\bigtriangleup$
\end{Remark}

The possible values of the index of a Kato surface with one branch
only depend on its black nodes.
It will be shown in the following section that
the upper bound is the index of the sublattice spanned by the rational
curves in $H^2(S,\z)$, see \cite[3.14]{dl11}:

\begin{Proposition} \label{indeces}
Given any finite sequence of positive integers ${k_0,k_1,...,k_{p-1}},$ choose
another positive integer $m$ and consider a Kato surface $S$ with
$\mathrm{DlS}=[s_{k_0}s_{k_1}...s_{k_{p-1}}r_m].$
Then, $\mathrm{index}(S)=1$ if and only if
$$m =\cp(k_0,k_1,...,k_{p-1}) \t  -  \cp(k_1,k_2,...,k_{p-1})$$
for some $\t\in\n$ which will then be its tip multiplicity.

For all other $m$ we have $2\leq \mathrm{index}(S) \leq \cp (k_0,k_1,...,k_{p-1})$
and all these values for the index are attained for some suitably chosen $m$.
\end{Proposition}

\section{More branches}

In order to deal with the general case of an intermediate Kato surface we now change notation.
For the simple Dloussky sequence $[s_{k_0}s_{k_1}...s_{k_{p-1}}r_m]$  we set:
\be
\cp := \cp(k_0,k_1,...,k_{p-1})
\ee
\be
\cq := \cq(k_1,k_2,...,k_{p-1},m) := \cp(k_1,k_2,...,k_{p-1}) + m .
\ee

As usual $r$ denotes the multiplicity of the unique root $R$.
Recall what we have computed after setting the tip multiplicity $a_1=:\t$
\be r =(\cp +1)\t -\cq +m \ee
\be c_0 = (\cp +1)\t -\cq +1 . \ee
where $c_0$ is the multiplicity of the first cycle black node $C_0$
and notice that the two formulas agree when $m=1$ because this is equivalent to $C_0=R$.
In fact, we know from Duality lemma \ref{dl} that $c_0=a_1 +1$, and it follows from
(\ref{remarkable}) that $c_0 = r - m +1$;
we then found the tip multiplicity $a_1=\t$ by solving the equation
\be \t = (\cp +1)\t -\cq  \quad \textrm{ from which } \quad  \t = \frac\cq\cp .
\ee

\bigskip\bigskip
Notice that $\cp$ and $\cq$ both depend on $p$ variables.
The variables of $\cp$ are the multiplicities of all black nodes
while $\cq$ is independent of the first black node and depends on $m$ additively.

\bigskip

Now, a general intermediate Kato surface $S$ will have $N$ branches and
we write its Dloussky sequence as
\be \mathrm{DlS} =[\mathrm{DlS_1} \cdots \mathrm{DlS_N}] \label{notsimple}\ee
where each $\mathrm{DlS_f}=[s_{k_{0f}}s_{k_{1f}}...s_{k_{(p_{f}-1)f}}r_{m_{f}}]$
is a simple sequence as before and will be called a simple component of $\mathrm{DlS}$.

\smallskip
We can then write down all multiplicities of $-K=c_1(S)$ by using the same
procedure as in the simple case, because adjunction formula is local,
in the sense that it only involves the nodes of $\Gamma$ having
a common edge with the given node.
Let again $\t:=a_{11}$ be the multiplicity of the first tip,
then $c_{01}=\t +1$ will be the multiplicity of the first black
node of the cycle and denote by $r[1]=(\cp_1 +1)\t -\cq_1 +m_1$ the
one of the first root $R_1$, using obvious notations.

Now that we reached the first root the following happens: suppose at first that
$R_1$ is a white node or equivalently $m_1\geq 2$,
by the same proof as before
the multiplicities will go down, one by one, along the chain of $-2$ curves
in the regular sequence $r_{m_{1}}$ until they reach the next
black node in the cycle, namely $C_{02}\in\mathrm{Dl}(S_2)$.
Its multiplicity will then be
$c_{02}=(\cp_1 +1)\t -\cq_1 +1$ and we get that the second
tip multiplicity must be
\be a_{12}=(\cp_1 +1)\t -\cq_1 . \ee
Notice that these formulas hold unchanged even in the case $m_1=1$
or equivalently $R_1=C_{02}$.
Continuing in the same way we will get that
\be a_{13}=(\cp_2 +1)a_{12} -\cq_2 \ee
and so on until the cycle closes up and we get back
to the first tip $a_{11}=\t$ which will give us the
following linear equation for $\t$, showing in particular that $\t >0$ :

\be [\prod_{f=1}^N(\cp_f+1) -1]\t =
 \sum_{h=1}^N [\cq_h \prod_{f=h+1}^N(\cp_f+1)]. \label{osaka} \ee

We can now collect some consequences of this formula and its proof;
they show that the results obtained in the simple case generalize
to arbitrary intermediate Kato surfaces.

\begin{Theorem} Let $S$ be an intermediate Kato surface with
Dloussky sequence (\ref{notsimple}). Then
$$\mathrm{index}(S)= \frac{\prod_{f=1}^N(\cp_f +1) -1}
{g.c.d.\{ [\prod_{f=1}^N(\cp_f+1) -1 ],
\sum_{h=1}^N [\cq_h \prod_{f=h+1}^N(\cp_f+1)] \} }$$
\end{Theorem}


Notice that the right hand side of  (\ref{osaka}) depends upon a choice of first tip
while the $\mathrm{index}(S)$ does not.
The reader can easily check this fact algebraically.

\medskip
By the same argument as in the simple case we also
get restrictions on the number of branch components
of an intermediate Kato surface of index $1$.
Let $b_2(B)=\sum k_{(even)}f$ be the total number
of irreducible components of all the $N$ branches.
The number of cycle components equals
$b_2(C)=\sum k_{(odd)f} + m $ where we set

\begin{Definition} \rm
Let $S$ be an intermediate Kato surface with $N\geq 1$ branches and
Dloussky sequence (\ref{notsimple}), with simple components
Dl$S_f=[s_{k_{0f}}s_{k_{1f}}...s_{k_{(p_{f}-1)f}}r_{m_{f}}]$.
It will be important to consider the sum of the lengths of all regular sequences
$$m:=\sum_{f=1}^{N} m_f .$$
\end{Definition}

The following result extends \ref{shortbranch} to the general case.
It says in particular that the number of branch components
cannot exceed the number of cycle components for an intermediate Kato surface
of index $1$.

\begin{Corollary} Let $S$ be an intermediate Kato surface of $\mathrm{index}(S)=1$.
Then,  $$b_2(B) \leq m.$$
 with equality   if and only if each of the simple components Dl$S_f$ of Dl$S$
is of the form $[s_{m_{f}} r_{m_{f}}]$ or $[3s_{k_{f}}2]$.
In particular, $b_2(B) \leq b_2(C)$ with equality if and only if each
Dl$S_f=[s_{m_{f}} r_{m_{f}}]$.
\end{Corollary}

\it Proof. \rm
Assume first that there is a white root. After a cyclic permutation of the simple
components of DlS we can suppose that $R_1^2 =-2$ -- i.e. $m_1 \geq 2$.
This root will meet the first branch $B_1$ and the first piece of cycle
$C_{01}+C_{11}+\cdots +C_{{\beta_1}{1}}$.
As in the simple case the branch multiplicities form a
strictly increasing sequence of positive integers, from
$\t$ up to the root multiplicity, denoted by $r[1]$.
After that, in the cyclic order of the cycle, the multiplicities
will go down $1$ by $1$ along a chain of $(-2)$--curves
until the first black node $C_{02}$ of the second simple component DlS$_2$
and will start to increase up to its root $R_2$, see Corollary \ref{slopes}.

The situation will go on in the same way but when we reach a black root,
call it $R_{l}$ -- i.e. $m_{l}=1$ and $R_{l}=C_{0l}$ -- a new phenomenon occurs.
First of all, its self-intersection number $R^2_l=-(k_{0l} +2)$ can be arbitrarily negative;
furthermore, setting as usual $g_{1l}:=c_{1l}-c_{0l}$,
we easily see by adjunction that the slope $g_{1l}$ satisfies the initial conditions of
Theorem \ref{chase}:
$g_{1l}=k_{0l}(c_{0l} -1) -1=k_{0l}a_{0l} -1$ which is certainly non-negative since $a_{0l}\in N$.
This shows  that, when the index is $1$, the multiplicity $c_{0l}$ of a black root $C_{0l}$
is not a local extremum for the cycle multiplicities.

The conclusion is that the absolute minimum of the cycle multiplicities
occurs at a black node $C_{0i}$ which is not a root while
the maximal multiplicity occurs at a white root.
In between these two values, the multiplicities are
non-decreasing along chains  of type
$C_{0f}+C_{1f}+\cdots +C_{{\beta_f}{f}}+R_f$
and are strictly decreasing as slowly as possible
along chains of type
$R_f+C_{{(\beta_f +2)}{f}}+\cdots +C_{0(f+1)}$;   each of
these chains has length $(m_f -1)$ so that some of them maybe empty
and this happens precisely when $R_f$ is a black root.

This argument shows that the total decrease of the cycle multiplicities
is $m-N\geq 0$ and has to coincide with the total growth,
because the cycle closes up:
$$m-N= (r[1] -c_{01}) +(r[2] -c_{02}) +\cdots +(r[N] -c_{0N}).$$
Finally, by Theorem \ref{chase},  $a_{1f} +1=c_{0f}$ for all $f=1,\dots,N$
so that the total growth along the branches is precisely $m$.

The result now follows from the simple observation that, because the index is $1$,
the branch multiplicities form a strictly increasing sequence of integers
by Corollary \ref{slopes} (in particular the branch slopes are all $\geq 1$)
and therefore the number of branch nodes is at most $m$: $b_2(B)\leq m$.

The same argument as in the simple case  shows that
 equality holds if and only if each simple component
is of the form $[s_kr_k]$ or else $[3s_k2]$.
As $m\leq b_2(C)$ always holds we conclude that index$(S)=1$ implies
$b_2(B) \leq b_2(C)$ with equality only when each simple component is
of the form $[s_kr_k]$.
\hfill    $\Box$
\medskip


\section{Applications}

Our work is related and motivated by the work of several authors. Let
$$  \mathrm{DlS}=[\mathrm{DlS_1} \dots \mathrm{DlS_N}]  $$
denote a Dloussky sequence of a Kato surface $S$ with $N$ branches.
It determines the directed dual graph of $S$ and conversely
two Kato surfaces have the same directed dual graph $\Gamma$ if and only if their
Dloussky sequences differ by a cyclic permutation of simple components
$\mathrm{Dl}S_f$.
We will then denote by $\mathcal{C}_\Gamma$ the set of Kato surfaces
with the same directed dual graph $\Gamma$.

\medskip
The index of a Kato surface only depends on $\Gamma$
and is invariant by unramified coverings.
Some motivations for studying it are the following:
first of all by \cite{do99}, $\mathrm{index}(S)=1$
is a necessary condition for an intermediate Kato surface $S$
to admit a twisted holomorphic vector field or a twisted anticanonical section.
By \cite[prop.2]{agg} it is also a necessary condition
for $S$ to admit bi-Hermitian metrics.
Our results then give precise obstructions on $\Gamma$ 
for the existence of these holomorphic sections or metric structures.

Although we don't know of any example of bi-Hermitian metrics
on intermediate Kato surfaces, it follows from \cite{do99}
that if $\mathrm{index}(\Gamma)=1$ then there are $S\in \mathcal{C}_\Gamma$
with a holomorphic vector field; as well as different $S\in \mathcal{C}_\Gamma$
with a holomorphic anticanonical section.

\bigskip
We now pass to present applications of our results starting from a relation with a
recent work of Dloussky \cite{dl11} in which the author computes the discriminant
$k=k(S)$ of the singularity obtained by contracting
to a point the maximal curve $D$ of a Kato surface $S$.
Our first aim is to indicate how $k$ is related to $\mathrm{index}(S)$.
We start with the following lemma which shows that our multilinear forms
$\cp(X_1,...,X_n)$ coincide with the simplest version of Dloussky
polynomials, denoted by $\cp_{\{n\}}(X_1,...,X_n)$ and defined as follows

\begin{Definition} \cite[p.35]{dl11} \rm Let $X_1,...,X_n$ denote a set of variables
and define the following polynomial 
$$\cp_{\{n\}}(X_1,...,X_n) := \sum \prod_{j\notin B} X_j$$
where $B\subseteq \{0,1,...,n\}$ ranges over all possible subsets
(including $\emptyset$) which can be written as \it disjoint \rm union of the
following building blocks: $\{n\}; \{1,2\}; \{2,3\};...; \{n-1,n\}$.
\end{Definition}

\begin{Lemma} For every $n\in\n$, the multilinear form $\cp(X_1,...,X_n)$
inductively defined in (\ref{p}) coincides with the Dloussky polynomial $\cp_{\{n\}}(X_1,...,X_n)$.
\end{Lemma}
\it Proof. \rm The first step of induction is easily verified. Then,
we write the Dloussky polynomial as sum
of monomials containing $X_n$, plus the rest:
$\cp_{\{n\}}(X_1,...,X_n)=
 X_n[\sum \prod_{j\notin B'} X_j] + \sum \prod_{j\notin B''} X_j$.
Here $B=\{n\}\cup B''$, and therefore $B''\subseteq \{0,1,...,n-1\}$ can be
written as disjoint union of $\{n\}\setminus \{n\}=\emptyset;
 \{1,2\}; \{2,3\};...; \{n-1,n\}\setminus \{n\} = \{n-1\}.$
We conclude that the rest - i.e. all the monomials without $X_n$ -
are just $\cp_{\{n-1\}}(X_1,...,X_{n-1})$ which, by induction
equals $\cp(X_1,...,X_{n-1})$.

Now we come to monomials containing $X_n$:
they are all of the form $X_n\prod_{j\notin B'} X_j$
where $B'\subseteq \{1,2,,...,n-1\}$ can be written as disjoint union of
$\emptyset;  \{1,2\}; \{2,3\};...; \{n-2,n-1\}.$
Notice that there is no singleton in these building blocks.
To complete the proof it only remains to show that our inductively
defined polynomials $f$ in (\ref{f}) actually coincide with

$$\tilde{f}(X_1,...,X_{n-1}):=
  \left\{ \begin{array}{ll}
      \sum \prod_{j\notin B'} X_j     & \mbox{if $n$ is even};\\  & \\
       \sum \prod_{j\notin B'} X_j    +1  & \mbox{ otherwise }.\end{array} \right. $$

This can be done by induction. The first step is easily verified;
assume by induction that $\tilde{f}(X_1,...,X_k)=f(X_1,...,X_k)$ for all $k<n-1$.
To finish the proof, we write $\tilde{f}$ as sum of monomials which contain the
variable $X_{n-1}$ and those which do not.

This corresponds to the following two possibilities for $B'$, we start
form the latter one:

(i) $n-1 \in B'$ in which case $\{n-2,n-1\} \subset B'$
   therefore the building blocks are $\{1,2\};...;\{n-3,n-4\}$,
   by induction their contribution is precisely $f(X_1...,X_{n-3})$.

(ii) $n-1 \notin B'$ and each monomial is of the form
    $X_{n-1}\prod_{j\notin B'}X_j$ where $B'\subseteq \{0,1,...,n-1\}$
    is disjoint union of  $\emptyset;  \{1,2\}; \{2,3\};...; \{n-2,n-1\}$
    which equals $X_{n-1}f(x_1,...,X_2)$, again by induction.

The proof is complete because we have shown that
$$\tilde{f}(X_1,...,X_{n-1})=X_{n-1}f(X_1,...,X_{n-2}) + f(X_1,...,X_3)$$
 \hfill  $\Box$  \bigskip

From the above lemma and the results of \cite{dl11}
we immediately get the following

\begin{Corollary} Let $S$ be a Kato surface with $N$ branches and set
$k:=\prod_{f=1}^N(\cp_f+1)$ -- in particular $k=\cp +1$ when $S$ is simple -- then,

(i) $k-1$ is the index of the sublattice spanned by the rational
    curves in $H^2(S,\z)$.

(ii) $(k-1)^2$ is the determinant of the intersection
     matrix of the rational curves.

(iii) Each factor $\cp_f +1$ equals the determinant of the
  intersection matrix of the corresponding branch $B_f$.

(iv) $k$ is the twisting coefficient of the singularity obtained by contracting
  the maximal curve $D\subset S$ to a point.

(v) $k$ is also the exponent ``$k$" of the contracting germ of $S$.

\end{Corollary}
\bigskip

In order to effectively compute the relevant exponents in the contracting germ of $S$
we prove:

\begin{Lemma}\label{division} For every $n\in\n$, our multilinear forms (\ref{p}) satisfy
the following identity:
$$\cp(X_1,...,X_n) = X_1[\cp(X_2,...,X_n) +1] + \cp(X_3,...,X_n).$$
\end{Lemma}
\it Proof. \rm By the previous result
$\cp(X_1,...,X_n) = \sum \prod_{j\notin B} X_j$ with the same $B$ as before.
We separate the monomials containing $X_1$ and write
$\cp=X_1[\sum\prod_{j\notin D'}X_j] + \sum\prod_{j\notin D''}X_j$.
Where $1\in D''$ which actually implies $\{1,2\}\subset D''$,
therefore $\sum\prod_{j\notin D''}X_j = \cp(X_3,...,X_n)$.

Next, we look at $X_1[\sum\prod_{j\notin D'}X_j]$ and notice that
it will always contain the monomial $X_1$ alone, whether $n$ is even or odd.
We see from this that $\sum\prod_{j\notin D'}X_j= 1+\cp(X_2,...,X_n)$
because $1\notin D'$ says that $D'$ can be written as disjoint union of $\{n\},\{2,3\},...,\{n-1,n\}$.
\hfill  $\Box$ \bigskip

We are now ready to explain some relations among our results
on $\mathrm{index}(S)$, the contracting germ of $S$ and the dimension of the logarithmic
moduli space of the pair $(S,D)$ where $D\subset S$ is the maximal curve.

By the work of Dloussky \cite{dl84} the complex structure of a Kato surface $S$
is completely determined by a contracting polynomial germ of $\c^2$
which can be put in the following normal form \cite{fa00}, \cite[4.2]{ot08}
\be \phi(z,\xi) = (\lambda \xi^s z + P(\xi) +c\xi^{\frac{sk}{k-1}}, \xi^k) \label{favre} \ee
\be \mathrm{ where \quad} P(\xi) = \xi^j +c_{j+1}\xi^{j+1}+...+c_s \xi^s  \label{poly} \ee

The germ is said to be in \it pure \rm normal form when the coefficient $c=0$,
this condition is automatically satisfied if $\mathrm{index}(S)>1$ \cite{ot08}.

It is shown in \cite{ot08} how the contracting germ determines
the Dloussky sequence $\mathrm{DlS}$ \cite[sect. 6]{ot08} and this is
used there  to describe the logarithmic moduli space of the pair $(S,D)$
because varying the coefficients $(\lambda, c_{j+1},\dots ,c_s,c)$ corresponds
to changing the complex structure of $S$, leaving the maximal curve $D$ fixed.
When $\mathrm{DlS}=[s_{k_{0}}\dots s_{k_{(p-1)}} r_m]$ is simple,
every exponent between $j$ and $s$ appears in the contracting germ
and therefore the number of non-zero coefficients is $s-j+1=m$ in the pure case,
and is $m+1$ otherwise.  Moreover, we can now show

\begin{Proposition} Let $\mathrm{DlS}=[s_{k_{0}}\dots s_{k_{(p-1)}} r_m]$
be a simple Dloussky sequence.
Then, for the relevant exponents of its contracting germ we have the following:
$k=\cp(k_0,...,k_{p-1}) +1$, while $j=\cp(k_1,...,k_{p-1}) +1$ and $s=\cq=j+m-1$.
Furthermore, $\mathrm{g.c.d.}(j,k)=1$.
\end{Proposition}
\it Proof. \rm The result on $k$ has been already established.
To compute $j$ and $s$ we use the algorithm of \cite[p.336]{ot08}
where, in their notation, $\alpha_i = k_{i-1}$, for $i=1,...,p$.

We know that $k=\cp(k_0,...,k_{p-1}) +1$ therefore,
by lemma \ref{division}, we get $j=\cp(k_1,...,k_{p-1})+1$
and $\beta_1=\cp(k_2,...,k_{p-1})+1$.
In general, $s=j+m-1$ \cite[p.337]{ot08} so that $s=\cq$.

Finally, we check $\mathrm{g.c.d.}(j,k)=1$, see \cite{ot08}.
By repeated applications of \ref{division} we have:
$\mathrm{g.c.d.}(j,k)=\cdots =\mathrm{g.c.d.}(\cp(k_0,k_1)+1, k_1+1)=\mathrm{g.c.d.}(k_0k_1+k_0+1,k_1+1)=1$.
\hfill $\Box$

\begin{Corollary}
The following relations hold on any Kato surface with one branch:

The tip multiplicity $\t=\cq/\cp =s/(k-1)$.

The root multiplicity $r=(\cp +1)\t -\cq +m = k\t-s+m = \frac{sk}{k-1} +1-j$.


The exponent $s=(k-1)\t $.

The exponent $\frac{sk}{k-1}=k\t$.

\end{Corollary}

When $S$ has more branches, say $N>1$ branches,
it is shown in \cite[Section 5]{ot08} that the associated germ
is the composition of \it simple \rm germs $\phi_i$ with $i=1,\dots,N$;
meaning that $\mathrm{gcd}(j,k)=1$ for each simple component.
Furthermore the germ $\phi$ of $S$ is pure if and only if each $\phi_i$ is pure.
We then see from the proof of \cite[Proposition 5.10]{ot08} that
$k=\prod_{i=1}^N k_i$ while the exponents $s$ and $j$ are not uniquely defined
but depend on the composition order. For example if $N=2$ then
$j=j_1k_2$ and $s=s_1k_2+s_2$.
From this discussion and our formula (\ref{osaka}) it then easily follows:

\begin{Proposition} Let $\mathrm{DlS} =[\mathrm{Dl}S_1\dots\mathrm{Dl}S_N]$
be a Dloussky sequence of a Kato surface $S$ with contracting germ $\phi$ then,
$\t=\frac{s}{k-1}$; $\mathrm{index}(S)=\frac{k-1}{\mathrm{gcd}(k-1,s)}$ as in
the simple case, while $\mathrm{gcd}(k,j)>1$ when $N>1$.
Furthermore, if  DlS$_f = [s_{k_{0f}},\ldots ,s_{k_{(p_f-1)f}}r_{m_f}]$
indicate the simple components of Dl$S$ for $f=1,\dots ,N$ then,
the number of non-zero coefficients in $\phi$ equals
$m(=\sum_{i=1}^N m_i)$ when $\phi$ is pure and is $m+1$ otherwise.
\end{Proposition}
\medskip

We conclude with a different approach for computing the number of moduli of the pair.

Let $S$ be an intermediate Kato surface with maximal curve $D$
with $b$ irreducible components $D_i, 1\le i\le b$,
where $b=b_2(S)$ is the second Betti number of $S$.
Suppose that $S$ is associated with the Dloussky sequence
DlS = [DlS$_1 \ldots $ DlS$_N$]
where each DlS$_f$ is of the form  $[s_{k_{0f}} \ldots s_{k_{(p_f-1)f}}r_{m_f}]$
and $N$ is the number of branches of $D$.
In the following proposition we compute directly
the dimension of the tangent space $H^1(S,\Theta_S(-\log D))$
of the Kuranishi space of deformations
of the pair $(S,D)$ in terms of the Dloussky sequence, or equivalently,
the directed dual graph of $S$.
We recall that
$\delta = \delta_S:= \dim H^0(S,\Theta_S)\le 1$ and $\delta = 1$ if and only if
$S$ admits a non-trivial $\c$-action, which in turn occurs only when
index$(S)=1$ (but not vice versa) (cf.\ \cite{do99}).

\begin{Proposition}\label{}
Let $S$ be as above. Then we have
$\dim H^1(S,\Theta_S(-\log D))= m + \delta$.
\end{Proposition}

{\em Proof}.
We write $h^1(F)=h^1(X,F)$ for any sheaf on a complex manifold  $X$.
By Nakamura \cite[\S 3]{na90} we have $h^2(\Theta_S)=h^2(\Theta_S(-\log D))=0$
and by \cite{do99}         \newline
$h^0(\Theta_S)=h^0(\Theta_S(-\log D))=\delta$.
Hence by Riemann-Roch formula we have $h^1(\Theta_S) = 2b + \delta$.
Now we consider the sheaf exact sequence
\begin{equation}\label{i}
 0 \ra \Theta_S(-\log D) \ra \Theta_S  \ra \oplus_{l=1}^b N_l \ra 0,    \\
\end{equation}
where $N_l = N_{D_l/S}$ is the normal bundle of $D_l$ in $S$.
Suppose first that all $D_l$ are smooth.
Then since $n_l:= \deg N_l$ is negative, we have $h^0(N_l)=0$ for all $l$,
where deg denotes the degree. Hence the associated cohomology exact sequence
yields the following short exact sequence:
\begin{equation}\label{ir}
0 \ra H^1(\Theta_S(-\log D)) \ra H^1(\Theta_S)
\ra \oplus_{l=1}^b H^1(N_l) \ra 0.
\end{equation}
Therefore applying Serre duality to each $H^1(N_l)$ we have
\begin{equation}\label{irr}
h^1(\Theta_S(-\log D)) = 2b+\delta - \sum_{l=1}^b (-n_l-1).
\end{equation}
Now we write the index set $L:=\{1,\ldots ,b\}$ as the disjoint union $L=I\amalg J$,
where $n_l=-2$ for $l\in I$ and $n_l\le -3$ for $l\in J$. In other words,
$D_l,l\in L$, corresponds
to a white (resp.\ black) node if $l\in I$ (resp.\ $J$).
For $l\in J$ we set $k_l = n_l+2$. Then we get
\begin{equation}\label{}
\sum_l(-n_l-1)
= \sum_{l\in I}(2-1) + \sum_{l\in J}(k_l+1) = \# I + \sum_{l\in J} (k_l+1),
\end{equation}
where $\# I$ is the number of elements of $I$.
Now from the relation of Dloussky sequence and the weighted graph
associated to  $D$ explained in \cite[\S 6]{ot08} (cf.\ also the figure on p.4 of this paper) we have the following relation:   \newline
1) $\sum_{l\in J} k_l + \sum_{f=1}^N m_f = b$ and 2) $\#I = \sum_{l\in J} (k_l-1) + \sum_{f=1}^N m_f$.
In particular we have
\begin{eqnarray*}
\sum_l(-n_l-1)
& = &  \sum_{l\in J}\{(k_l-1)+(k_l+1)\} + \sum_{f=1}^N m_f \\
& = &  2 (\sum_{l\in J}{k_l} + \sum_{f=1}^N m_f)
- \sum_{f=1}^N m_f = 2b - \sum_{f=1}^N m_f.
\end{eqnarray*}
Hence we finally get
$h^1(\Theta_S(-\log D)) = 2b +\delta - (2b - m) = m + \delta$.

It remains to consider the case where $D$ contains a singular irreducible curve $D$.
In this case the unique ``cycle'' on $S$ is irreducible and
is identified with a rational curve, say  $D_1$, with a single rational double point; it has a unique branch consisting of smooth irreducible components
with self-intersection number $-2$; its Dloussky (simple) sequence is
given by $[s_{b-1},r_1]$.

Let $u: \tilde{S} \ra  S$ be the unique unramified double covering.
Since $S$ is an intermediate Kato surface, any irreducible components of the
maximal curve $\tilde{D}$ on $\tilde{S}$ is a smooth rational curve.
Now by the Riemann-Roch theorem the alternating sum $\chi(S,\Theta_S(-\log D)):= \sum_{i=0}^2 (-1)^ih^i(S,\Theta_S(-log D))$ is given by the Todd characteristic $T(S,\Theta_S(-\log D))$ of
the pair $(S,\Theta_S(-\log D))$ and the same is true for
$\chi(\tilde{S},\Theta_{\tilde{S}}(-\log \tilde{D}))$.
Since $u$ is unramified, we have
$\Theta_{\tilde{S}}(-\log \tilde{D}))\cong u^*\Theta_S(-\log D))$ and hence
$T(\tilde{S},\Theta_{\tilde{S}}(-\log D))=2u^*T(S,\Theta_S(-\log D))$.
Thus we get
$\chi(\tilde{S},\Theta_{\tilde{S}}(-\log \tilde{D})) = 2\chi(S,\Theta_S(-\log D))$.
Taking into account of the vanishing of $h^2$-terms on both sides
as noted above, we get   \newline
$h^1(\tilde{S},\Theta_{\tilde{S}}(-\log \tilde{D}))-\tilde{\delta}
= 2(h^1(S,\Theta_S(-log D))-\delta)$ with $\tilde{\delta}$ denoting $\delta$ on $\tilde{S}$.

On the other hand, $u^{-1}(D_1)$ is a cycle of two rational curves $\tilde{D}_{1s}, s=1,2$, with self-intersection number $(\tilde{D}_{1s})^2 = D_1^2-2\le -3$.
From this we conclude that $\tilde{S}$ has $2N$ branches and $2\sum_{f=1}^N m_f
=\sum_{\tilde{f}=1}^{2N} \tilde{m}_{\tilde{f}}$. (Actually, $N=1, m_f =\tilde{m}_{\tilde{f}}=1$.)
Together with the first part of the proof this gives the desired equality
also in this exceptional case.
\hfill  $\Box$
\medskip

\begin{Remark}  \rm
The fact that the dimension of the logarithmic moduli space equals the total length of
the regular sequences is geometrically explained as follows.
In the correspondence between the blowing-up sequence and the Dloussky sequence
each entry of singular sequence corresponds to the blowing up with center one of
the nodes of the previously produced exceptional curves and hence has no moduli,
while for the entry of regular sequence the blowing up occurs on the general points
and hence each contributes to one dimensional moduli.
\hfill  $\bigtriangleup$    \end{Remark}

We can now compute the dimension of the tangent space $H^1(S, \Theta_S(-D))$
of the Kuranishi family of deformations of $S$
which are trivial along the maximal curve $D$.

\begin{Proposition} Let $S$ be a Kato surface, then
$\dim H^1(S, \Theta_S(-D)) = b+m-\epsilon$. Where $b$ is the second Betti number of $S$,
$m$ is the total number of components in the regular sequences of $S$ while $\epsilon=1$
when $S$ admits both a holomorphic vector field and a tip $A_{1h}$ of multiplicity $a_{1h}=1$;
otherwise $\epsilon=0$.
\end{Proposition}
\it Proof. \rm
We consider the exact sequence of sheaves in $S$
$$0\to \Theta_S(-D) \to \Theta_S(-\log D) \to T \to 0$$
where $T$ is the tangent sheaf of the maximal curve $D$.
Since $D$ is the union of $b$ irreducible components $D_i$ meeting transversally and
which we can assume to be smooth rational curves,
we see that $T=\oplus_{i=1}^b T_i$.
Now, $T_i\cong \mathcal{O}_{\mathbb{CP}_1}$ when $D_i$ meets two other components,
$T_i\cong \mathcal{O}_{\mathbb{CP}_1}(1)$ if $D_i$ is a tip and the only other possibility
is that $D_i$ is a root in which case $T_i\cong \mathcal{O}_{\mathbb{CP}_1}(-1)$.

It easily follows form this and from the previous proposition that the Kuranishi
family is unobstructed with $\quad \dim H^1(S, \Theta_S(-D)) = b+m \quad$
whenever the following map is isomorphic:
$H^0(S, \Theta_S(-D)) \to H^0(S,\Theta_S(-\log D))$.

The only other possibility is that $S$ has a holomorphic vector field, automatically
tangent to $D$ along $D$, which does not vanish identically on $D$.
In this case we have to subtract $\epsilon=1$ dimensions and notice that index$(S)=1$.
Finally, it is shown in \cite{do99} that a holomorphic vector field $\theta$
vanishes exactly on the divisor $D_{\theta} := -K-D$.
Therefore $\theta$ is not identically zero on the maximal curve $D$ if and only if $1$
is the minimal multiplicity of $-K$. But we have already shown that the absolute
minimum for the multiplicities of a surface of index $1$ occurs at some of its tips.
\hfill $\Box$



\newcommand{\bysame}{\leavevmode\hbox to3em{\hrulefill}\,}

\bigskip
\noindent
A. Fujiki -- Research Institute for Mathematical Sciences, Kyoto University. \\
M. Pontecorvo --  Dipartimento di Matematica e Fisica, Roma Tre University


\begin{thebibliography}{99}

\bibitem[AGG99]{agg} Apostolov, V., Gauduchon, P., Grantcharov, G.,
    Bi-Hermitian structures on complex surfaces.
    {\em Proc. London Math. Soc.} (3) {\bf 79} (1999), no. 2, 414--428.

\bibitem[Dl84]{dl84} Dloussky, G.,
    Structure des surfaces de Kato.
    \it M\'emoires dela S.M.F. \rm   {\bf 112} 14 (1984).

\bibitem[Dl06]{dl06} Dloussky, G.,
    On surfaces of class  VII$_0$ with numerically
    anticanonical divisor. {\em Amer. J. Math.} {\bf 128} (2006), no. 3, 639--670.

\bibitem[Dl11]{dl11} Dloussky, G.
  Quadratic forms and singularities of genus one or two.
  {\it Ann. Fac. Sci. Toulouse} {\bf XX} (2011) 15 -- 69.

\bibitem[DO99]{do99} Dloussky, G. ; Oeljeklaus, K.
  Vector fields and foliations an compact surfaces of class $VII_0$.
  {\it Ann. Inst. Fourier} {\bf 49} (1999) 1503 -- 1545.

\bibitem[Fa00]{fa00} Favre, C.
  Classification of 2-dimensional contracting rigid germs.
  { \it J. Math. Pure Appl.} {\bf 79} 475 -- 514 (2000).


\bibitem[Ka77]{ka77} Kato, Ma.
   Compact complex manifolds containing global spherical shells.
    {\it I. Proceedings of the International Symposium on Algebraic Geometry}
    (Kyoto Univ., Kyoto, 1977), pp. 45–84, Kinokuniya Book Store, Tokyo, 1978.

\bibitem[Na89]{na89} Nakamura, I.
 Towards classification of non-K\"ahlerian couplex surfaces.
 {\it Sugaku Expositions} {\bf 2} (1989) 209 -- 229.

\bibitem[Na90]{na90} \bysame,
    On surfaces of class VII$_0$ with curves, II.
    {\em T\^ohoku Math. J.}, {\bf 42} (1990) 475-516.

\bibitem[OT08]{ot08}  Oeljeklaus, K.; Toma, M.
    Logarithmic moduli spaces for surfaces of class $VII$.
     {\it Math. Ann.} {\bf 341} (2008) 323 -- 345.
     arXiv:math/0701840v2 [math.CV] (2009).

\bibitem[Te10]{te10} Teleman, A.
Istantons and holomorphic curves on class VII surfaces.
{\it Ann. Math.} {\bf 172} (2010) 1749 -- 1804.

\bibitem[Te13]{te13}  \bysame.
A variation formula for the determinant line bundle.
Compact subspaces of moduli spaces of stable bundles over class VII surfaces.
{\it 	arXiv:1309.0350}.

\end{thebibliography}
\end{document}